\input amstex
\magnification=\magstep1 
\baselineskip=13pt
\documentstyle{amsppt}
\vsize=8.7truein \CenteredTagsOnSplits \NoRunningHeads
\def\UU{\Cal U}
\def\cut{\operatorname{cut}}
 \topmatter
 
\title Computing the partition function for graph homomorphisms \endtitle 
\author Alexander Barvinok and Pablo Sober\'on \endauthor
\address Department of Mathematics, University of Michigan, Ann Arbor,
MI 48109-1043, USA \endaddress
\email barvinok$\@$umich.edu, psoberon$\@$umich.edu  \endemail
\date May 2015 \enddate
\thanks  The research of the first author was partially supported by NSF Grants DMS 0856640 and DMS 1361541.
\endthanks 
\abstract 
We introduce the partition function of edge-colored graph homomorphisms, of which the usual partition function of graph homomorphisms is a specialization, and present an efficient algorithm to approximate it in a certain domain. Corollaries include efficient algorithms for computing weighted sums approximating the number of $k$-colorings and the number of independent sets in a graph, as well as an efficient procedure to distinguish pairs of edge-colored graphs with many color-preserving homomorphisms $G \longrightarrow H$ from pairs of graphs that need to be substantially modified to acquire a color-preserving homomorphism $G \longrightarrow H$.
\endabstract

\endtopmatter

\document

\head 1. Introduction and main results \endhead

\subhead (1.1) Graph homomorphism partition function  \endsubhead Let $G=(V, E)$ be an undirected graph with set $V$ of vertices and set 
$E$ of edges, without multiple edges or loops, and let $A=\left(a_{ij}\right)$ be a $k \times k$ symmetric complex matrix. The
{\it graph homomorphism partition function} is defined by 
$$P_G(A)=\sum_{\phi: V \rightarrow \{1, \ldots, k\}} \prod_{\{u, v\} \in E} a_{\phi(u) \phi(v)}. \tag1.1.1$$
Here the sum is taken over all maps $\phi: V \longrightarrow \{1, \ldots, k\}$ and the product is taken over all edges in $G$.

The function $P_G(A)$ encodes many interesting properties of the graph $G$, and, not surprisingly, is provably hard to compute except in a few special cases, see \cite{C+13} and references therein. For example, if $A$ is the adjacency matrix of  an undirected graph $H$ with vertices $1, \ldots, k$, that is, if
$$a_{ij}=\cases 1 &\text{if \ } \{i, j\} \text{\ is an edge of \ } H \\ 0 &\text{otherwise,} \endcases$$
then $P_G(A)$ is the number of {\it homomorphisms} of $G$ into $H$, that is, the number of maps 
$\phi: V \longrightarrow \{1, \ldots, k\}$ such that $\{\phi(u), \phi(v)\}$ is an edge of $H$ whenever $\{u, v\}$ is an edge of $G$.

Here are some examples of a particularly interesting choices of the matrix $A$, see also Section 5.3 of \cite{Lo12} for more.

\subsubhead (1.1.2) Colorings \endsubsubhead
If the $k \times k$ matrix $A$ is defined by
$$a_{ij}=\cases 1 &\text{if \ } i \ne j \\ 0 &\text{if\ } i=j,\endcases$$
then $P_G(A)$ is the number of $k$-{\it colorings} of $G$, that is, the number of ways to color the vertices of $G$ into $k$ colors so that the endpoints of every edge of $G$ have different colors. Indeed, each $k$-coloring of $G$ contributes 1 to 
$P_G(A)$ in (1.1.1) via the map $\phi: V \longrightarrow \{1, \ldots, k\}$ that maps the vertices colored in the $i$-th color into $i$.
The smallest $k$ for which a $k$-coloring of $G$ exists is called the {\it chromatic number} of $G$. Approximating the chromatic number of a given graph within a factor $|V|^{1-\epsilon}$ is NP-hard for any fixed $\epsilon >0$ \cite{FK98},  \cite{Zu07}.

\subsubhead (1.1.3) Independent sets \endsubsubhead
Suppose that $k=2$ and that $A$ is defined by 
$$a_{ij}=\cases 0 &\text{if \ } i=j=1 \\ 1 &\text{otherwise.} \endcases$$
Then $P_G(A)$ is the number of {\it independent sets} in $G$, that is, the number of subsets $U \subset V$ of vertices such that no two vertices of $U$ span an edge of $G$. Indeed, each independent set $U$ contributes 1 to $P_G(A)$ in (1.1.1) via the map $\phi: V \longrightarrow \{1, 2\}$ such that $\phi^{-1}(1)=U$.

\subsubhead (1.1.4) Maximum cut \endsubsubhead 
Suppose that $k=2$. For $0 < \epsilon < 1$, let us define $A=A_{\epsilon}$ by 
$$a_{ij}= \cases \epsilon &\text{if\ } i=j \\ 1 &\text{if\ } i \ne j \endcases$$
and let us consider the value of $\epsilon^{-|E|}P_G\left(A_{\epsilon}\right)$. Every map $\phi: V \longrightarrow \{1, 2\}$ in (1.1.1) is uniquely defined by the subset $S \subset V$ such that $S=\phi^{-1}(1)$. For a subset $S \subset V$ we define the {\it cut}
associated with $S$ by 
$$\cut_G(S)=\left| \{u, v\} \in E: \ u \in S,\ v \notin S \right|.$$
Then 
$$\epsilon^{-|E|}P_G\left(A_{\epsilon}\right) = \sum_{S:\ S \subset V} \epsilon^{-\cut_G(S)}.$$
Let 
$$\mu(G)=\max_{S:\ S \subset V} \cut_G(S)$$
be the maximum cut associated with a subset $S$ of vertices.
Then 
$$\epsilon^{-\mu(G)} \ \leq \ \epsilon^{-|E|} P_G\left(A_{\epsilon}\right) \ \leq \ 2^{|V|} \epsilon^{-\mu(G)}$$
and hence
$${\ln P_G\left(A_{\epsilon}\right) \over \ln (1/\epsilon)} + |E| - |V| {\ln 2 \over \ln (1/\epsilon)}\ \leq \ \mu(G)\ \leq \ {\ln P_G\left(A_{\epsilon}\right) \over \ln (1/\epsilon)} + |E|.$$
In particular, computing $P_G\left(A_{\epsilon}\right)$ for a sufficiently small, yet fixed, $\epsilon >0$, we can approximate $\mu(G)$ within an additive error of $\delta |V|$ for an arbitrarily small $\delta >0$, fixed in advance.

\subhead (1.2) Partition function of edge-colored graph homomorphisms \endsubhead Let $G=(V, E)$ be a graph as above and let $B=\left(b_{ij}^{uv}\right)$ be a $|E| \times {k(k+1) \over 2}$ complex matrix with entries indexed by edges $\{u, v\} \in E$ and unordered pairs $1 \leq i, j \leq k$. Technically, we should have written $b^{\{u, v\}}_{\{i, j\}}$, but we write just 
$b^{uv}_{ij}$, assuming that 
$$b^{uv}_{ij}=b^{vu}_{ij}=b^{vu}_{ji}=b^{uv}_{ji}.$$
We define the {\it edge-colored graph homomorphism partition function} by
$$Q_G(B)=\sum_{\phi: V \rightarrow \{1, \ldots, k\}} \prod_{\{u, v\} \in E} b^{uv}_{\phi(u) \phi(v)}, \tag1.2.1$$
where, as in (1.1.1), the sum is taken over all maps $\phi: V \longrightarrow \{1, \ldots, k\}$ and the product is taken over all edges of $G$. If $A=\left(a_{ij}\right)$ is a $k \times k$ symmetric matrix and we define $B$ by 
$$b^{uv}_{ij}=a_{ij} \quad \text{for all} \quad \{u, v\} \in E,$$
then 
$$Q_G(B)=P_G(A),$$
so $P_G$ defined by (1.1.1) is a specialization of $Q_G$ defined by (1.2.1). 

Let $H$ be an undirected simple graph with $k$ vertices and suppose that the edges of $G$ and $H$ are colored.
Let us define
$$b^{uv}_{ij}=\cases 1&\text{if\ } \{u, v\} \ \text{and} \ \{i, j\} \ \text{are edges of the same color} \\ &\text {of\ } G \ \text{and} \ H \  
\text{respectively} \\ 0 &\text{otherwise.}\endcases$$
Then $Q_G(B)$ is the number of \text{\it edge-colored homomorphisms} of $G$ into $H$, that is, the number of maps 
$\phi: V \longrightarrow \{1,\ldots, k\}$ such that for every edge $\{u, v\}$ of $G$, the pair $\{\phi(u), \phi(v)\}$ is an edge of $H$ of the same color, cf., for example, \cite{AM98}.

\subhead (1.3) Our results \endsubhead 
Let $\Delta(G)$ denote the largest degree of a vertex of $G$ (the degree of a vertex of a graph is the number of edges incident to the vertex). We present a deterministic algorithm, which, given a graph $G=(V, E)$, an $\epsilon > 0$ and a (real or complex) 
$|E| \times {k(k+1) \over 2}$ matrix $B=\left(b_{ij}^{uv}\right)$
 such that 
$$\left| 1 - b_{ij}^{uv} \right| \ \leq \ {\gamma \over \Delta(G)} \quad \text{for all} \quad \{u,v\} \in E \quad \text{and}
\quad 1 \leq  i,j \leq k,$$
where $\gamma >0$ is an absolute constant,
computes the value of $Q_G(B)$ within relative error $\epsilon$ in $\bigl(|E| k\bigr)^{O(\ln |E|-\ln \epsilon)}$ time. 
We can choose $\gamma=0.34$, if $\Delta(G) \geq 3$ we can choose $\gamma=0.45$, and for all sufficiently large $\Delta(G)$ we can choose $\gamma=0.54$. 

Consequently, we obtain an algorithm of $\bigl(|E| k\bigr)^{O(\ln |E|-\ln \epsilon)}$ complexity to approximate $P_G(A)$ for any 
$k \times k$ symmetric matrix $A=\left(a_{ij}\right)$ which satisfies
$$\left| 1- a_{ij} \right| \ \leq \ {\gamma \over \Delta(G)} \quad \text{for all} \quad 1\leq i, j \leq k.$$
This allows us to compute efficiently various ``soft" relaxations of ``hard" combinatorial quantities of interest. Here are 
the corresponding modification of Examples 1.1.2 and 1.1.3.

In Example 1.1.2, let us define the $k \times k$ matrix $A$ by 
$$a_{ij}=\cases 1+\gamma/\Delta(G) &\text{if\ } i \ne j \\ 1-\gamma/\Delta(G) &\text{if\ } i=j. \endcases$$
Then the value of 
$$\left(1 + {\gamma \over \Delta(G)}\right)^{-|E|} P_G(A) \tag1.3.1$$
represents the weighted  sum over all $k^{|V|}$ possible colorings of the vertices of $G$ into $k$ colors, where each proper coloring is counted with weight 1, whereas a coloring for which $w$ edges are miscolored (that is, have their endpoints colored with the same color) is counted with weight
$$\left(1+{\gamma \over \Delta(G)}\right)^{-w} \left(1-{\gamma \over \Delta(G)}\right)^w \ \leq \ \exp\left\{-{2 \gamma w \over \Delta(G)}\right\}. \tag1.3.2$$

In Example 1.1.3, let us define the $2 \times 2$ matrix $A$ by
$$a_{ij}=\cases 1-\gamma/\Delta(G) &\text{if\ } i=j=1 \\ 1+\gamma/\Delta(G) &\text{otherwise.} \endcases$$
Then the value of (1.3.1) represents the weighted sum over all $2^{|V|}$ subsets of vertices of the graph $G$, where each independent set is counted with weight 1, whereas a set whose vertices span $w$ edges of $G$ is counted with weight (1.3.2).

Let us restrict ourselves to the class of graphs of bounded degree, with $\Delta(G) \leq 3$, say. Then our result implies that the value of the partition function $P_G(A)$ can be efficiently approximated as long as $1-\delta \leq a_{ij} \leq 1+\delta$ for all $i$ and $j$, where 
$0< \delta <1 $ is an absolute constant (we can choose $\delta=0.11$). It is tempting to conjecture that for any $0< \delta < 1$, fixed in advance, the value of $P_G(A)$ can be efficiently approximated. This, however, cannot be so unless NP-hard problems can be solved by a quasi-polynomial algorithm. Indeed, approximating the maximum cut in $G$ satisfying $\Delta(G) \leq 3$ within a certain absolute constant factor $\beta_0 >1$ is known to be NP-hard \cite{BK99}. The problem remains NP-hard if we further restrict ourselves to connected graphs satisfying $\Delta(G) \leq 3$. In this case the maximum cut is at least $|V|-1$ and the construction of Section 1.1.4 shows that for some fixed 
$\epsilon >0$ approximating $P_G\left(A_{\epsilon}\right)$ within some fixed factor $\beta_1>1$ is an NP-hard problem.

We note that for any positive $A$ the problem of computing $P_G(A)$ exactly is $\#P$-hard unless $\operatorname{rank} A=1$, in which case the problem admits a polynomial time algorithm \cite{BG05}.

Computing $Q_G(B)$ allows us to distinguish pairs of edge-colored graphs with many color-preserving homomorphisms 
$G \longrightarrow H$ from pairs which are sufficiently far from having a single color-preserving homomorphism.
Indeed, given edge-colored graphs $G$ and $H$, let us define $B=\left(b^{uv}_{ij}\right)$ by
$$b^{uv}_{ij}=\cases 1 +{\gamma \over \Delta(G)} &\text{if\ } \{u, v\} \ \text{and} \ \{i, j\} \ \text{are edges of the same color} \\ &\text {of\ } G \ \text{and} \ H \  
\text{respectively} \\ 1-{\gamma \over \Delta(G)} &\text{otherwise.}\endcases$$
Then the value of 
$$\left(1+{\gamma \over \Delta(G)}\right)^{-|E|} Q_G(B) \tag1.3.3$$
represents the weighted sum over all $k^{|V|}$ maps $\phi: V \longrightarrow \{1, \ldots, k\}$, where each color-preserving homomorphism is counted with weight 1 and a map $\phi$ which does not map some $w$ edges of $G$ onto the identically colored edges of $H$ is counted with weight (1.3.2) at most. 

Let us choose some positive integer $w$.
Hence if every map $\phi$ does not map some $w$ edges of $G$ onto the identically colored edges of $H$, the value of (1.3.3) does not exceed 
$ k^{|V|} e^{-2 \gamma w / \Delta(G)}$ (in this case, we say that $G$ and $H$ are sufficiently far from having a color-preserving 
homomorphism $G \longrightarrow H$).
If, however, the probability that a random map $\phi$ is a color-preserving homomorphism is at least
$2e^{-2 \gamma w/\Delta(G)}$, then the sum (1.3.3) is at least 
$2 k^{|V|} e^{-2 \gamma w/\Delta(G)}$ (in this case we say that there are sufficiently many color-preserving homomorphisms).
Computing the value of $Q_G(B)$ within relative error $0.1$, say, we can tell apart these two cases. The most interesting situation is when $G$ is almost regular, so $|E| \approx 0.5 |V| \Delta(G)$ and $w \approx \epsilon |E|$ for some fixed 
$\epsilon >0$, in which case ``many" may still mean that the probability to hit a color-preserving homomorphism at random is exponentially small.

\subhead (1.4) The idea of the algorithm \endsubhead Let $J$ denote the $|E| \times {k (k+1) \over 2}$ matrix filled with 1s.
Given a $|E| \times {k(k+1) \over 2}$ matrix $B=\left(b^{uv}_{ij}\right)$, where $\{u, v\} \in E$ and 
$1 \leq i, j \leq k$, we consider the univariate function 
$$f(t)=\ln Q_G\bigl(J+t(B-J)\bigr),\tag1.4.1$$
so that 
$$f(0)=\ln Q_G(J)=|V| \ln k \quad \text{and} \quad f(1)=\ln Q_G(B).$$
Hence our goal is to approximate $f(1)$ and we do it by using the Taylor polynomial expansion of $f$ at $t=0$:
$$f(1) \approx f(0)+ \sum_{m=1}^n {1 \over m!} {d^m \over dt^m} f(t) \Big|_{t=0}. \tag1.4.2$$
It turns out that the approximation (1.4.2) can be computed in $(|E| k)^{O(n)}$ time. We present the algorithm in Section 2.
The quality of approximation (1.4.2) depends on the location of {\it complex} zeros of $Q_G$.

\proclaim{(1.5) Lemma} Suppose that there is a real $\beta >1$ such that 
$$Q_G\bigl(J+z(B-J)\bigr) \ne 0 \quad \text{for all} \quad z \in {\Bbb C} \quad \text{satisfying} \quad |z| \leq \beta.$$
Then the right hand side of (1.4.2) approximates $f(1)$ within an additive error of 
$${|E| \over (n+1)\beta^n\left(\beta-1\right)}.$$
\endproclaim
In particular, for a fixed $\beta >1$, to ensure an additive error of $0 < \epsilon < 1$, we can choose 
$n=O\left(\ln |E| -\ln \epsilon\right)$, which would result in the algorithm for approximating $Q_G(B)$ within 
relative error $\epsilon$ in $(|E| k)^{O(\ln |E| -\ln \epsilon)}$ time. We prove Lemma 1.5 in Section 2.

It remains to identify a class of matrices $B$ for which the number $\beta >1$ of Lemma 1.5 exists. We prove the following result.

\proclaim{(1.6) Theorem} There exists an absolute constant $\alpha >0$ such that for any undirected graph $G$ and any complex $|E| \times {k(k+1)\over 2}$ matrix $B=\left(b^{uv}_{ij}\right)$ satisfying 
$$\left| 1-b_{ij}^{uv} \right| \ \leq \ {\alpha \over \Delta(G)} \quad \text{for all} \quad \{u, v\} \in E \quad \text{and} \quad
1 \leq  i,j \leq k,$$
where $\Delta(G)$ is the largest degree of a vertex of $G$, one has 
$$Q_G(B) \ne 0.$$
One can choose $\alpha=0.35$, if $\Delta(G) \geq 3$ one can choose $\alpha=0.46$ and if $\Delta(G)$ is sufficiently large, one can choose $\alpha=0.55$.
\endproclaim

We prove Theorem 1.6 in Section 3. Theorem 1.6 implies that if 
$$\left| 1- b_{ij}^{uv}\right| \ \leq \ {0.34 \over \Delta(G)} \quad \text{for all} \quad \{u, v\} \in E \quad \text{and} \quad 
1 \leq  i, j \leq k,$$
we can choose $\beta=35/34$ in Lemma 1.5 and hence obtain an algorithm which computes $Q_G(B)$ within relative error 
$\epsilon$ in $(|E| k)^{O(\ln |E| -\ln \epsilon)}$ time. Similarly, if $\Delta(G) \geq 3$ and 
$$\left| 1- b_{ij}^{uv} \right| \ \leq \ {0.45 \over \Delta(G)} \quad \text{for all} \quad \{u, v\} \in E \quad \text{and}  \quad 
1 \leq i, j \leq k,$$
we can choose $\beta=46/45$ and if 
$$\left| 1- b_{ij}^{uv}\right| \ \leq \ {0.54 \over \Delta(G)} \quad \text{for all} \quad \{u, v\} \in E \quad \text{and} \quad
1 \leq i, j \leq k,$$
and $\Delta(G)$ is sufficiently large, (namely, if $\Delta(G) \geq 30$) we can choose $\beta=55/54$.

A similar approach was used earlier to compute the permanent of a matrix \cite{B15a} and the partition function for cliques of a given size in a graph \cite{B15b}. While the algorithm of Section 2 and Lemma 1.5 are pretty straightforward modifications of the corresponding results of \cite{B15a} and \cite{B15b}, the proof of Theorem 1.6 required new ideas.

\head 2. The algorithm \endhead

\subhead (2.1) The algorithm for approximating the partition function \endsubhead We present an algorithm, which, 
given a $|E| \times {k(k+1)\over 2}$ matrix $B=\left(b_{ij}^{uv}\right)$, computes the approximation (1.4.2) for the function $f$ defined by (1.4.1).
Let
$$g(t)=Q_G\bigl(J + t(B-J)\bigr), \tag2.1.1$$
so $f(t)=\ln g(t)$. Hence
$$f'(t)={g'(t) \over g(t)} \quad \text{and} \quad g'(t)=g(t) f'(t).$$
Therefore, for $m \geq 1$, we have 
$${d^m \over dt^m} g(t)\Big|_{t=0}=\sum_{j=0}^{m-1} {m-1 \choose j} \left({d^j \over dt^j} g(t)\Big|_{t=0}\right)\left({d^{m-j} \over dt^{m-j}} f(t)\Big|_{t=0} \right) \tag2.1.2$$
(we agree that the $0$-th derivative of $g$ is $g$).
We note that $g(0)=k^{|V|}$. If we compute the values of 
$${d^m \over dt^m} g(t)\Big|_{t=0} \quad \text{for} \quad m=1, \ldots, n, \tag2.1.3$$
then the formulas (2.1.2) for $m=1, \ldots, n$ provide a non-degenerate triangular system of linear equations that allows us to compute
$${d^m \over dt^m} f(t) \Big|_{t=0} \quad \text{for} \quad m=1, \ldots, n.$$
Hence our goal is to compute the values (2.1.3).
We have 
$${d^m \over dt^m} g(t)\Big|_{t=0}=\sum_{\phi: V \rightarrow \{1, \ldots, k\}} \sum \Sb I=\bigl(\{u_1, v_1\}, \\ \qquad \ldots, \\ \qquad \{u_m, v_m\}\bigr) \endSb \left(b_{\phi(u_1) \phi(v_1)}^{u_1 v_1}-1 \right) \ldots \left(b_{\phi(u_m) \phi(v_m)}^{u_m v_m}-1\right),$$
where the inner sum is taken over all ordered sets $I$ of $m$ distinct edges $\{u_1, v_1\}$, $\ldots$, $\{u_m, v_m\}$ of $G$. Let  $S(I)$ be the set of all distinct vertices among $u_1, v_1$, $\ldots$, $u_m, v_m$. Then 
$$\split &{d^m \over dt^m} g(t)\Big|_{t=0}\\ &\quad=\sum_I  k^{|V|-|S(I)|}   \sum_{\phi:\ S(I) \rightarrow \{1, \ldots, k\}} 
\left(b_{\phi(u_1) \phi(v_1)}^{u_1 v_1}-1 \right) \ldots \left(b_{\phi(u_m) \phi(v_m)}^{u_m v_m}-1\right), \endsplit$$
where the outer sum is taken over not more than $|E|^m$ ordered sets $I$ of $m$ distinct
edges $\{u_1, v_1\}, \ldots, \{u_m, v_m\}$ of $G$ and the inner sum is taken over not more than $k^{2m}$ maps 
$\phi: S(I) \longrightarrow \{1, \ldots, k\}$. Hence the complexity of computing the approximation (1.4.2) is 
$\left(|E| k\right)^{O(n)}$ as claimed.

\subhead (2.2) Proof of Lemma 1.5 \endsubhead The function $g(t)$ defined by (2.1.1) is a polynomial of degree $d \leq |E|$ and $g(0)=k^{|V|} \ne 0$, so we factor
$$g(z)=g(0) \prod_{i=1}^d \left(1-{z \over \alpha_i}\right),$$
where $\alpha_1, \ldots, \alpha_d \in {\Bbb C}$ are the roots of $g(z)$. By the condition of Lemma 1.5, we have 
$$\left| \alpha_i \right| \ \geq \ \beta \ > \ 1 \quad \text{for} \quad i=1, \ldots, d.$$ Therefore,
$$f(z)=\ln g(z) =\ln g(0) + \sum_{i=1}^d \ln \left(1 -{z \over \alpha_i}\right) \quad \text{for} \quad |z| \leq 1, \tag2.2.1$$
where we choose the branch of $\ln g(z)$ that is real at $z=0$. Using the standard Taylor expansion, we obtain
$$\ln \left(1- {1 \over \alpha_i}\right) = -\sum_{m=1}^n {1 \over m} \left({1 \over \alpha_i}\right)^m + \zeta_n,$$
where 
$$\left| \zeta_n\right|= \left| \sum_{m=n+1}^{+\infty} {1 \over m} \left({1 \over \alpha_i}\right)^m \right| \ \leq \ {1 \over (n+1) \beta^n (\beta-1)}.$$
Therefore, from (2.2.1) we obtain
$$f(1) =f(0) +\sum_{m=1}^n \left(-{1 \over m} \sum_{i=1}^d \left({1 \over \alpha_i}\right)^m \right) + \eta_n,$$
where
$$\left| \eta_n \right| \ \leq \ {|E| \over (n+1)\beta^n (\beta-1)}.$$
It remains to notice that 
$$-{1\over m} \sum_{i=1}^d \left({1 \over \alpha_i}\right)^m ={1 \over m!} {d^m \over dt^m} f(t) \Big|_{t=0}.$$
{\hfill \hfill \hfill} \qed

\head 3. Proof of Theorem 1.6 \endhead

For a $0 < \delta < 1$, we define the polydisc $\UU(\delta) \subset {\Bbb C}^{k(k+1)|E|/2}$ by 
$$\UU(\delta)=\Bigl\{Z=\left(z^{uv}_{ij}\right): \quad \left| 1-z^{uv}_{ij}\right| \leq \delta \quad \text{for all} \quad \{u, v\} \in E \quad 
\text{and} \quad 1 \leq i, j \leq k \Bigr\}.$$ Thus we have to prove that for 
$\delta=\alpha/\Delta(G)$, where $\alpha >0$ is an absolute constant, we have $Q_G(Z) \ne 0$ for all $Z \in \UU(\delta)$.

\subhead (3.1) Recursion \endsubhead For a sequence of distinct vertices $W=\left(v_1, \ldots, v_m\right)$ of the graph $G$ and a sequence 
$L=\left(l_1, \ldots, l_m\right)$ of not necessarily distinct numbers $1 \leq l_1, \ldots, l_m \leq k$, we define 
$$Q^W_L(Z) = \sum \Sb \phi: V \rightarrow \{1, \ldots, k\} \\ \phi\left(v_1\right)=l_1, \ldots, \phi\left(v_m\right)=l_m \endSb 
\prod_{\{u, v\} \in E} z^{uv}_{\phi(u) \phi(v)}$$
(we suppress the graph $G$ in the notation).
In words: we restrict the sum (1.2.1) defining $Q_G(Z)$ onto the maps $\phi: V \longrightarrow \{1, \ldots, k\}$ that map selected vertices $v_1, \ldots, v_m$ of $G$ into preassigned indices $l_1, \ldots, l_m$. We denote $|W|$ the number of vertices in $W$ and $|L|$ the number of indices in $L$ (hence we have $|W|=|L|$).

We denote by $(W, u)$ a sequence $W$ appended by $u$ (distinct from all previous vertices in $W$) and by $(L, l)$ a sequence $L$ appended by $l$ (not necessarily distinct from all previous indices in $L$).
Then for any sequence $W$ of distinct vertices, for any $u$ distinct from all vertices in $W$ and for any sequence $L$ of indices such that $|L|=|W|$, we have 
$$Q^W_L(Z)=\sum_{l=1}^k Q^{(W, u)}_{(L, l)}(Z). \tag3.1.1$$
When $W$ and $L$ are both empty, then $Q^W_L(Z)=Q_G(Z)$.

We start with a geometric inequality.
\proclaim{(3.2) Lemma} Let $x_1, \ldots, x_n \in {\Bbb R}^2$ be non-zero vectors such that for some 
$0 \leq \alpha < 2\pi/3$ the angle between any two vectors $x_i$ and $x_j$ does not exceed $\alpha$. Let 
$x=x_1 + \ldots + x_n$. Then 
$$\|x\| \ \geq \ \left(\cos {\alpha \over 2} \right) \sum_{i=1}^n \|x_i\|.$$
\endproclaim
\demo{Proof} We note that $0$ is not in the convex hull of any three vectors $x_i, x_j, x_k$, since otherwise the angle between some two of those three vectors would have been at least $2 \pi/3$. The Carath\'eodory Theorem implies that $0$ is not in the convex hull of $x_1, \ldots, x_n$ and hence the vectors lie in an angle of at most $\alpha$ with vertex at the origin. 
Let us consider the bisector of that angle and the orthogonal projections of each $x_i$ onto the bisector. The length of the orthogonal projection of each $x_i$ is at least $\|x_i\| \cos(\alpha/2)$ and hence the length of the orthogonal projection of
$x_1 + \ldots + x_n$ is at least $\left(\|x_1\| + \ldots + \|x_n\|\right) \cos(\alpha/2)$. Since the vector $x_1 + \ldots + x_n$ is at least as long as its orthogonal projection, the proof follows.
{\hfill \hfill \hfill} \qed
\enddemo

Lemma 3.2 was suggested by Boris Bukh \cite{Bu15}. It replaces a weaker bound of 
$\sqrt{\cos \alpha}\left(\|x_1\| + \ldots + \|x_n\|\right)$, assuming that $\alpha \leq \pi/2$, of an earlier version of the paper.

Our proof of Theorem 1.6 is based on the following two lemmas.

\proclaim{(3.3) Lemma} Let $\tau >0$ be real, let $W$ be a sequence of distinct vertices of $G$, let 
$u$ be a vertex distinct from the vertices in $W$ and let $L$ be a sequence of not necessarily distinct numbers from
the set $\{1, \ldots, k\}$ such that $|L|=|W|$.  Suppose that for all 
$Z \in \UU(\delta)$ and for all $1 \leq l \leq k$, we have  
$$Q^{(W, u)}_{(L, l)} (Z) \ne 0$$
and, moreover,  
$$\left|Q^{(W, u)}_{(L, l)}(Z) \right| \  \geq \ {\tau \over \Delta(G)}  \sum \Sb v:\ \{u, v\} \in E  \\ j:\ 1 \leq j \leq k \endSb
\left| z^{u v}_{l j} \right| \left| {\partial \over \partial z^{uv}_{l j}}Q^{(W, u)}_{(L, l)}(Z) \right|.$$
Then, for any two $1 \leq l, m \leq k$ and any $A \in \UU(\delta)$, the angle between two complex numbers 
$$Q^{(W, u)}_{(L, l)}(A) \quad \text{and} \quad Q^{(W, u)}_{(L, m)}(A),$$
interpreted as vectors in ${\Bbb R}^2={\Bbb C}$, does not exceed 
$$\theta ={2 \delta \Delta(G)  \over \tau (1-\delta)}.$$
\endproclaim
\demo{Proof} Since $Q^{(W, u)}_{(L, l)}( Z) \ne 0$ for all $Z \in \UU(\delta)$, we can and will consider a branch of 
$\ln Q^{(W, u)}_{(L, l)}( Z)$ for $Z \in \UU(\delta)$. Then 
$${\partial \over \partial z^{uv}_{l j}}\ln Q^{(W, u)}_{(L, l)}( Z) ={\partial \over \partial z^{uv}_{l j}}Q^{(W, u)}_{(L, l)}(Z)\Big/
Q^{(W, u)}_{(L, l)}(Z)$$ and since 
$$\left| z^{xy}_{ij} \right| \ \geq \ 1-\delta \quad \text{for all} \quad x, y, i, j$$
we conclude that 
$$\sum\Sb v:\ \{u, v\} \in E  \\ j:\ 1 \leq j \leq k \endSb \left| {\partial \over \partial z^{uv}_{l j}} \ln  Q^{(W, u)}_{(L, l)}(Z) \right| 
\ \leq \ {\Delta(G) \over \tau(1-\delta)} \quad \text{for all} \quad Z \in \UU(\delta).$$
Given $A \in \UU(\delta)$, $A=\left(a^{xy}_{ij}\right)$, and $1 \leq l , m \leq k$, we define $B \in \UU(\delta)$, $B=\left(b^{xy}_{ij}\right)$,
by 
$$b^{uv}_{lj}=a^{uv}_{mj} \quad \text{for all} \quad v \in V \quad \text{such that} \quad \{u, v\} \in E \quad \text{and all} \quad
1 \leq j \leq k$$
and 
$$b^{xy}_{ij}=a^{xy}_{ij} \quad \text{in all other cases}.$$
Then 
$$Q^{(W, u)}_{(L, l)}(B)=Q^{(W, u)}_{(L, m)}(A)$$ 
and hence 
$$\split
& \left| \ln Q^{(W, u)}_{(L, l)}(A) -  \ln Q^{(W, u)}_{(L, m)}(A) \right| = \left| \ln Q^{(W, u)}_{(L, l)}(A) -  \ln Q^{(W, u)}_{(L, l)}(B) \right|  \\ &\quad \leq 
\max_{Z \in \UU(\delta)} \sum\Sb v:\ \{u, v\}\in E  \\ j:\ 1 \leq j \leq k \endSb \left| {\partial \over \partial z^{uv}_{l j}} \ln  Q^{(W, u)}_{(L, l)}(Z) \right| \times \max\Sb v \in V:\ \{u, v\} \in E \\ j: \ 1 \leq j \leq k \endSb
 \left| a^{uv}_{lj}-b^{uv}_{lj}\right|   \ \leq \ {2\delta \Delta(G) \over \tau (1-\delta)}, \endsplit $$
where the last inequality follows since $\left| a^{xy}_{ij}-b^{xy}_{ij}\right| \leq 2\delta$ for all $A, B \in \UU(\delta)$.
The proof now follows.
{\hfill \hfill \hfill} \qed
\enddemo

\proclaim{(3.4) Lemma} Let $0 < \theta < 2\pi/3$ be a real number, let $W$ be a sequence of distinct vertices and let $L$ be a sequence of not necessarily distinct indices from the set $\{1, \ldots, k\}$ such that $|L|=|W|$. Suppose that for any 
$Z \in \UU(\delta)$, for every $v \in V$ distinct from the vertices of $W$, and for every $1 \leq i, j \leq k$ we have 
$$Q^{(W, v)}_{(L, i)}(Z),\quad Q^{(W, v)}_{(L, j)}(Z) \ne 0$$ and that the angle between 
$$Q^{(W, v)}_{(L, i)}(Z) \quad \text{and} \quad  Q^{(W, v)}_{(L, j)}(Z),$$
considered as vectors in ${\Bbb R}^2 ={\Bbb C}$, does not exceed $\theta$.

Let $W=\left(W', u\right)$ and $L=\left(L', l \right)$. Then for all $Z \in \UU(\delta)$ we have 
$$\left| Q^{W}_{L} (Z) \right| \ \geq \ {\tau \over \Delta(G)} \sum \Sb v: \ \{u, v\} \in E \\ j: \ 1 \leq j \leq k \endSb 
\left|z^{uv}_{lj}\right| \left| {\partial \over \partial z^{uv}_{lj}} Q^{W}_{L} (Z) \right|,$$
where
$$\tau = \cos {\theta \over 2}.$$
\endproclaim 
\demo{Proof} Let $v$ be a vertex of $G$ such that $\{u, v \} \in E$.
If $v$ is an element of $W'$ then 
$${\partial \over \partial z^{uv}_{lj}} Q^{W}_{L} (Z) ={1 \over z^{uv}_{lj}} Q^{W}_{L} (Z),$$
provided $j$ is the element in the $L'$ sequence which corresponds to $v$ and 
$${\partial \over \partial z^{uv}_{lj}} Q^{W}_{L} (Z) =0$$
if the element in the $L'$ sequence corresponding to $v$ is not $j$.

If $v$ is not an element of $W'$ then 
$${\partial \over \partial z^{uv}_{lj}} Q^{W}_{L} (Z)={\partial \over \partial z^{uv}_{lj}} Q^{(W, v)}_{(L, j)}
={1 \over z^{uv}_{lj}} Q^{(W, v)}_{(L, j)} (Z).$$
Denoting by $d_0$ the number of vertices $v$ in the sequence $W'$ such that $\{u, v\} \in E$, we obtain 
$$\sum \Sb v: \ \{u, v\} \in E \\ j: \ 1 \leq j \leq k \endSb 
\left|z^{uv}_{lj}\right| \left| {\partial \over \partial z^{uv}_{lj}} Q^{W}_{L} (Z) \right| = d_0 \left| Q^{W}_{L} (Z)\right| +
\sum \Sb v \text{\ not in \ } W' \\ \{u, v\} \in E \\ 1 \leq j \leq k \endSb \left| Q^{(W, v)}_{(L, j)} (Z) \right|. \tag3.4.1$$
On the other hand, from (3.1.1) and Lemma 3.2, we conclude that for each $v$ not in the sequence $W'$, we have 
$$\left| Q^W_L(Z)\right| \ \geq \ \left(\cos {\theta \over 2}\right) \sum_{j=1}^k \left|Q^{(W, v)}_{(L, j)}(Z) \right|. \tag3.4.2$$
Denoting by $d_1$ the number of vertices $v$ not in the sequence $W'$ such that $\{u, v\} \in E$, we deduce from (3.4.1) and (3.4.2) that 
$$\sum \Sb v: \ \{u, v\} \in E \\ j: \ 1 \leq j \leq k \endSb 
\left|z^{uv}_{lj}\right| \left| {\partial \over \partial z^{uv}_{lj}} Q^{W}_{L} (Z) \right| 
\ \leq \ d_0 \left| Q^{W}_{L} (Z)\right| + {d_1 \over \cos {\theta \over 2}} \left| Q^W_L(Z)\right|,$$
from which the proof follows.
{\hfill \hfill \hfill} \qed
\enddemo

\subhead (3.5) Proof of Theorem 1.6 \endsubhead  One can see that for all sufficiently small $\alpha > 0$, the equation 
$$\theta = {2 \alpha \over (1-\alpha) \cos{\theta \over 2}}$$
has a solution $0 < \theta < 2 \pi/3$. Numerical computations show that one can choose 
$$\alpha =0.35 \quad \text{and} \quad \theta \approx 1.420166551.$$
Let 
$$\tau =\cos {\theta \over 2} \approx 0.7583075916.$$
Given a graph $G=(V, E)$ we define
$$\delta={\alpha \over \Delta(G)}$$
and prove by descending induction on $n=|V|, \ldots, 1$ the following three statements (3.5.1)--(3.5.3).
\bigskip
(3.5.1) For any sequence $W$ of $n$ distinct vertices of $G$, for every sequence $L$ of not necessarily distinct indices 
$1 \leq l \leq k$ such that $|W|=|L|$, for any $Z \in \UU(\delta)$, we have $Q^W_L(Z) \ne 0$;
\medskip
(3.5.2) Let $W$ be a sequence of $n$ distinct vertices of $G$ such that $W=(W', v)$ and let $L'$ be a sequence of not necessarily distinct indices $1 \leq l \leq k$ such that $|L'|=|W'|$. Then for every $1 \leq i, j \leq k$ and every $Z \in \UU(\delta)$, the angle between $Q^{(W', v)}_{(L', i)}(Z)$ and $Q^{(W', v)}_{(L', j)}(Z)$, interpreted as vectors in ${\Bbb R}^2={\Bbb C}$, does not exceed $\theta$;
\medskip
(3.5.3) Let $W$ be a sequence of $n$ distinct vertices of $G$ such that $W=(W', u)$ and let $L$ be a sequence of not necessarily 
distinct indices $1 \leq l \leq k$ such that $L=(L', l)$ and $|W|=|L|$. Then for all $Z \in \UU(\delta)$, we have
$$\left|Q^W_L(Z)\right| \ \geq \ {\tau \over \Delta(G)} \sum \Sb v:\ \{u, v\} \in E \\ j:\ 1 \leq j \leq k \endSb 
\left|z^{uv}_{lj} \right| \left| {\partial \over \partial z^{uv}_{lj}} Q^W_L(Z) \right|.$$
\bigskip
Suppose that $n=|V|$. If $W=\left(v_1, \ldots, v_n \right)$ and $L=\left(l_1, \ldots, l_n\right)$  then 
$$Q^W_L(Z) = \prod_{\left\{v_i, v_j\right\} \in E} z^{v_i v_j}_{l_i j_j} \ne 0,$$
so (3.5.1) holds. Moreover, denoting $\deg(v_n)$ the degree of $v_n$, we obtain
$$\sum \Sb v:\ \{v_n, v\} \in E \\ j:\ 1 \leq j \leq k \endSb 
\left|z^{v_n v}_{l_n j} \right| \left| {\partial \over \partial z^{v_n v}_{l_n j}} Q^W_L(Z)\right| = \deg\left(v_n\right) \left| Q^W_L(Z) \right|,$$
so (3.5.3) holds as well.

Statements (3.5.1) and (3.5.3) for sequences $W$ of length $n$ and Lemma 3.3 imply statement (3.5.2) for sequences 
$W$ and $L$ of length $n$.

Formula (3.1.1), Lemma 3.2 and statement (3.5.2) for sequences $W$ of length $n$ imply statement (3.5.1) for sequences $W$ of length $n-1$.

Statements (3.5.1) and (3.5.2) for sequences $W$ of length $n$ and Lemma 3.4 imply statement (3.5.3) for sequences  $W$ of length $n-1$.

This proves that (3.5.1)--(3.5.3) hold for sequences $W$ of length $1$. Formula (3.1.1), Lemma 3.2 and statement (3.5.2) for $n=1$ imply that $Q_G(Z) \ne 0$ for all $Z \in \UU(\delta)$.

We can improve the value of the constant $\alpha$ by defining $\theta$ as a solution to the equation
$$\theta ={2 \alpha \over \left(1-{\alpha \over \Delta(G)}\right) \cos {\theta \over 2}}.$$
Numerical computations show that one can choose $\alpha=0.55$ provided $\Delta(G) \geq 30$ and that one can choose $\alpha=0.46$ provided
$\Delta(G) \geq 3$.
\hfill \hfill \hfill \qed 

\head Acknowledgment \endhead

The authors are grateful to Boris Bukh for suggesting Lemma 3.2.


\Refs
\widestnumber\key{AAAA}

\ref\key{AM98}
\by N. Alon and T.H. Marshall
\paper Homomorphisms of edge-colored graphs and Coxeter groups
\jour Journal of Algebraic Combinatorics
\vol 8 
\yr 1998
\pages no. 1, 5--13
\endref

\ref\key{B15a}
\by A. Barvinok
\paper Computing the permanent of (some) complex matrices
\jour Foundations of Computational Mathematics
\paperinfo published online January 6, 2015, doi 10.1007/s10208-014-9243-7
\yr 2015
\endref

\ref\key{B15b}
\by A. Barvinok
\paper Computing the partition function for cliques in a graph
\jour Theory of Computing, to appear
\paperinfo  preprint \break {\tt  arXiv:1405.1974}
\yr 2015
\endref

\ref\key{BK99}
\by P. Berman and M. Karpinski
\paper On some tighter inapproximability results (extended abstract)
\bookinfo Automata, languages and programming (Prague, 1999)
\pages 200--209
\inbook Lecture Notes in Computer Science, 1644
\publ Springer
\publaddr Berlin
\yr 1999
\endref

\ref\key{Bu15}
\by B. Bukh
\paper Personal communication
\yr 2015
\endref

\ref\key{BG05}
\by A. Bulatov and M. Grohe
\paper The complexity of partition functions
\jour Theoretical Computer Science
\vol 348 
\yr 2005
\pages no. 2--3, 148--186
\endref

\ref\key{C+13}
\by J.-Y. Cai, X. Chen, and P. Lu
\paper Graph homomorphisms with complex values: a dichotomy theorem
\jour SIAM Journal on Computing 
\vol 42 
\yr 2013
\pages no. 3, 924--1029
\endref

\ref\key{FK98}
\by U. Feige and J. Kilian
\paper Zero knowledge and the chromatic number
\jour Journal of Computer and System Sciences
\vol 57 
\yr 1998
\pages  no. 2, 187--199
\endref

\ref\key{Lo12}
\by L. Lov\'asz
\book Large Networks and Graph Limits
\bookinfo American Mathematical Society Colloquium Publications, 60
\publ American Mathematical Society
\publaddr Providence, RI
\yr 2012
\endref

\ref\key{Zu07}
\by D. Zuckerman
\paper Linear degree extractors and the inapproximability of max clique and chromatic number
\jour Theory of Computing
\vol  3 
\yr 2007
\pages 103--128
\endref

\endRefs
\enddocument

\end